\input btxmac.tex

\long\def\comment#1{\relax}


\magnification=\magstep1 
\hsize=15truecm 
\vsize=22.3truecm
\hoffset=.45truecm    
\voffset=.3truecm


\hyphenation{Null-stel-len-satz Pro-po-si-tion}


\font\tenbb=msbm10     
\font\sevenbb=msbm7
\newfam\bbfam \def\bb{\fam\bbfam\tenbb}  
\textfont\bbfam=\tenbb
\scriptfont\bbfam=\sevenbb

\font\tenam=msam10     
\font\sevenam=msam7
\newfam\amfam   
\textfont\amfam=\tenam
\scriptfont\amfam=\sevenam

\font\twelverm=cmr12

\font\bn=cmb10

\font\titlefont=cmr12 at 13pt

\let\dinky=\sevenrm

\font\csc=cmcsc10            

\def\glq{\hbox{,\kern-.7pt,\kern.3pt}}  

\let\em=\it

\font\sevensl=cmsl8 at 7pt
\scriptfont\slfam=\sevensl

\let\subheadfont=\bf

\newcount\sectnum \sectnum=0
\newcount\itemnum

\def\numtag{\the\sectnum.\the\itemnum}

\newwrite\xrefs


\def\subhead#1{\removelastskip\goodbreak\bigskip\bigskip%
\advance\sectnum by 1 \itemnum=1
{\subheadfont  \the\sectnum.\enspace}{\subheadfont #1}\par%
\nobreak\medskip\noindent}

\def\nonumsubhead#1{\removelastskip\goodbreak\bigskip\bigskip%
\advance\sectnum by 1 \itemnum=1
{\subheadfont #1}\par%
\nobreak\medskip\noindent}

\let\tiny=\scriptscriptstyle


\thickmuskip=4.5mu plus 4mu minus 3mu  
\mathsurround=0.5pt

\pretolerance=200  \hyphenpenalty=1000 \exhyphenpenalty=5000

\baselineskip=14pt \parskip=1.5pt \parindent=1.5em
\dimen1=3pt \dimen2=2em 
 
\smallskipamount=4pt plus2pt minus1pt 
\medskipamount=7pt plus3pt minus2pt 
\bigskipamount=15pt plus5pt minus4pt 


\def\frac#1#2{{#1\over #2}}

\chardef\at="40  

\def\natn{{\bb N}}     
\def\intz{{\bb Z}}     
\def\ratq{{\bb Q}}     
\def\Dhat{{\hat D}} 
\def\Mhat{{\hat M}} 
\def\Ahat{{\hat A}} 
\def\an{a_1,\ldots,a_n}
\def\xn{x_1,\ldots,x_n}
\def\MnD{{\rm M}_n(D)}
\def\MnK{{\rm M}_n(K)}
\def\Mint{\IntD(\MnD)}
\def\MMint{\IntD[\MnD]}

\def\Pa{P_a}
\def\Mhata{(M\Ahat)_a}
\def\MAhata{(M\Ahat)_a}
\def\MAa{(MA)_a}

\def\quaternions{{\rm H}} 
\def\L{{\rm L}} 

\def\chiC{\chi_{\tiny C}}

\def\[{{\rm[}} \def\]{\/{\rm]}}  
\def\({{\rm(}} \def\){\/{\rm)}}  

\def\Int{\hbox{\rm Int}}  
\def\Intn{\hbox{\rm Int}^n}  
\def\IntD{\hbox{\rm Int}_D}  
\def\IntnD{\hbox{\rm Int}_D^n}  
\def\IntzeroD{\hbox{\rm Int}_D^0}  
\def\Im{\Int(A)}  
\def\Imn{\Intn(A)}  
\def\OK{{\cal O}_K} 


\let\tensor=\otimes 
\let\congr=\equiv    
\let\takeaway=\setminus  
\mathchardef\ideal="1945  
\mathchardef\eop="1903  

\newcount\codenumber \codenumber=`:
\loop\relax\if\codenumber>"FFF \advance\codenumber\by -"1000\repeat
\advance\codenumber by "3000
\mathcode`:=\codenumber

\def\label#1{
   \edef#1{\numtag}
   \write\xrefs{\def\string#1{#1}}
   }

\long\def\profess#1#2#3\endprofess {\ifdim\lastskip<\medskipamount%
  \removelastskip\medskip\penalty-500\fi
  \label#2
  \noindent{\bn\numtag\enspace#1.\ }{\sl#3\par}%
  \advance\itemnum by 1%
  \ifdim\lastskip<\smallskipamount\removelastskip\penalty-55\smallskip\fi}

\long\def\rmprofess#1#2#3\endrmprofess {\ifdim\lastskip<\medskipamount%
  \removelastskip\medskip\penalty-350\fi
  \label#2
  \noindent{\bn\numtag\enspace#1.\ }{\rm#3\par}%
  \advance\itemnum by 1%
  \ifdim\lastskip<\smallskipamount\removelastskip\penalty-100\smallskip\fi}

\def\proof{\noindent {\it Proof.\enspace}}
\def\endproof{\ $\eop$\par%
   \ifdim\lastskip<\medskipamount
   \removelastskip\penalty-50\medskip\fi}


\def\firstpagehead{\tenrm Appeared in J.~Algebra 373 (2013) 414-425\hfill \ }
\def\runninghead{\hfill 
{\dinky  INTEGER-VALUED POLYNOMIALS ON ALGEBRAS} \hfill }
\headline={\ifnum\pageno>1\runninghead \else\firstpagehead \fi}

\output={
          \shipout\vbox{\makeheadline
                        \vskip.1truein
                        \pagebody
                        \makefootline }
          \advancepageno
          \ifnum\outputpenalty>-20000 \else\dosupereject\fi
 }

\input mintxrefs.aux
\openout\xrefs=mintxrefs.aux

\centerline{\titlefont INTEGER\kern1.5pt\raise1.5pt\hbox{-}VALUED\ \ 
POLYNOMIALS\ \ ON\ \ ALGEBRAS}
\bigskip
\centerline{\twelverm Sophie Frisch}
\centerline{\tt frisch@TUGraz.at}
\bigskip\bigskip\bigskip

\hfuzz=1pt
\noindent{\csc Abstract.}
Let $D$ be a domain with quotient field $K$ and $A$ a $D$-algebra.
A polynomial with coefficients in $K$ that maps every element of $A$
to an element of $A$ is called integer-valued on $A$. For commutative
$A$ we also consider integer-valued polynomials in several variables.
For an arbitrary domain $D$ and $I$ an arbitrary ideal of $D$ we
show $I$-adic continuity of integer-valued polynomials on $A$.
For Noetherian one-dimensional $D$, we determine spectrum and Krull
dimension of the ring $\IntD(A)$ of integer-valued polynomials on $A$.
We do the same for the ring of polynomials with coefficients in
$M_n(K)$, the $K$-algebra of $n\times n$ matrices, that map every
matrix in $M_n(D)$ to a matrix in $M_n(D)$.
\medskip

\medskip

\subhead{Introduction}

Let $D$ be a domain with quotient field $K$ and $A$ a $D$-algebra,
such as, for instance, a group ring $D(G)$ or the matrix algebra
$M_n(D)$.

We are interested in the rings of polynomials
$$\IntD(A)=\{f\in K[x]\mid f(A)\subseteq A\},$$
and, if $A$ is commutative,
$$\IntnD(A)=\{f\in K[\xn]\mid f(A^n)\subseteq A\}.$$

Elements of the $D$-algebra $A$ are plugged into polynomials with
coefficients in $K$ via the canonical homomorphism 
$\iota_A:A\rightarrow K\tensor_D A$, $\iota_A(a)=1\tensor a$.

In the special case $A=D$ these rings are known as rings of
integer-valued polynomials, cf.~\cite{CaCh97ivp}. 
They provide natural examples of non-Noetherian Pr\"ufer rings
\cite{CCL01HDPivp,Lop98CIP}, and have been used for 
proving results on the $n$-generator property in Pr\"ufer rings
\cite{BoyKli06ngen}.
Also, integer-valued polynomials are useful for polynomial
interpolation of functions from $D$ to $D$ \cite{Fri99Iiv,CCF00id},
and satisfy other interesting algebraic conditions such as analogues
of Hilbert's Nullstellensatz \cite{CaCh97ivp,Fri01NSp}. 

These desirable properties of rings of integer-valued polynomials
have motivated the generalization to polynomials with coefficients
in $K$ acting on a $D$-algebra $A$ \cite{Fri05chh,LoWe12GRiv}. So far,
not much is known about rings of integer-valued polynomials on algebras.
We know that they behave somewhat like the classical rings of
integer-valued polynomials if the $D$-algebra $A$ is commutative. 
For instance, Loper and Werner \cite{LoWe12GRiv} have shown that 
$\Int_{\intz}(\OK)$ is Pr\"ufer.
If $A$ is non-commutative, however, the situation is
radically different. 
For instance, $\Int_{\intz}(M_2(\intz))$ is not Pr\"ufer \cite{LoWe12GRiv}, 
and is far from allowing interpolation \cite{Fri05chh}. 

We will describe the spectrum of $\IntD(A)$, for a one-dimensional
Noetherian ring $D$ and a finitely generated torsion-free $D$-algebra
$A$, in the hope that this will facilitate further research.
We will investigate more closely the special case of $A=\MnD$: we
determine a polynomially dense subset of $\MnD$ and describe
the image of a given matrix under the ring $\IntD(\MnD)$.

A different ring of integer-valued polynomials on the matrix 
algebra $M_n(D)$, consisting of polynomials with coefficients in
$M_n(K)$ that map matrices in $M_n(D)$ to matrices in $M_n(D)$,
has been introduced by Werner \cite{Wer12IVM}. We will show that
it is isomorphic to the algebra of $n\times n$ matrices over ``our''
ring $\IntD(\MnD)$ of integer-valued polynomials on $M_n(D)$ with
coefficients in $K$.

Before we give a precise definition of the kind of $D$-algebra $A$ for
which we will investigate $\IntD(A)$, a few examples. $D$ is always a
domain with quotient field $K$, and not a field. 

\rmprofess{Example}{\matrixexample}
For fixed $n\in\natn$, let $A=\MnD$ be the $D$-algebra of $n\times n$
matrices with entries in $D$ and
$$\IntD(\MnD)=\{f\in K[x]\mid \forall C\in \MnD:\; f(C)\in\MnD\}.$$
\endrmprofess

\rmprofess{Example}{\quaternionexample}
Let $\quaternions=\ratq + \ratq i + \ratq j + \ratq k$ be the
$\ratq$-algebra of rational quaternions, 
$\L= \intz + \intz i + \intz j + \intz k$ the $\intz$-subalgebra of
Lipschitz quaternions, and
$$\Int_\intz(\L)=\{f\in \ratq[x]\mid 
\forall z\in \L:\; f(z)\in\L\}.$$
\endrmprofess

\rmprofess{Example}{\groupringexample}
Let $G$ be a finite group, $K(G)$ and $D(G)$ the respective group rings,
and
$$\IntD(D(G))=\{f\in K[x]\mid \forall z\in D(G):\; f(z)\in D(G)\}.$$
If $G$ is commutative, we also consider
$$\IntnD(D(G))=\{f\in K[\xn]\mid \forall z\in D(G)^n:\; f(z)\in D(G)\}.$$
for $n\in\natn$, where $D(G)^n=D(G)\times\ldots\times D(G)$ denotes
the Cartesian product of $n$ copies of $D(G)$.
\endrmprofess

\rmprofess{Example}{\dedekindexample}
Let $D\subseteq A$ be Dedekind rings with quotient fields $K\subseteq F$,
and
$$\IntnD(A)=\{f\in K[\xn]\mid f(A^n)\subseteq A\}.$$
\endrmprofess

\rmprofess{Notation and conventions}{\generalconventions}
Throughout this paper, $D$ is a domain and not a field, $K$ the
quotient field of $D$, and $A$ a {\em torsion-free} $D$-algebra
{\em finitely generated} as a $D$-module. Since $A$ is faithful,
there is an isomorphic copy of $D$ embedded in $A$ by $d\mapsto d\/1_A$,
and we may assume $D\subseteq A$. 

Now let $B=K\tensor_D A$.
The natural homomorphisms $\iota_K:K\rightarrow K\tensor_D A$,
$\iota_K(k)=k\tensor 1$ and $\iota_A:A\rightarrow K\tensor_D A$,
$\iota_A(a)=1\tensor a$, are injective, since $A$ is a torsion-free
$D$ module.
We identify $K$ and $A$ with their isomorphic copies in $B$, which
allows us to evaluate polynomials with coefficients in $K$ at arguments
in $A$, and define
$$\IntD(A)=\{f\in K[x]\mid \forall a\in A:\; f(a)\in A\}$$
and for $n\in\natn_0$
$$\IntnD(A)=\{f\in K[\xn]\mid \forall \an\in A:\; f(\an)\in A\}.$$

To exclude pathological cases we require $K\cap A= D$.
\endrmprofess

\rmprofess{Remark}{\otherB}
Instead of $B=K\tensor_D A$, we could look at the canonically
isomorphic $A_{K\setminus\{0\}}$, the ring of fractions of $A$
with denominators in $K\setminus\{0\}$. The natural homomorphisms
$\iota_K:K\rightarrow A_{K\setminus\{0\}}$ and 
$\iota_A:A\rightarrow A_{K\setminus\{0\}}$ then take the form
$\iota_K({c\over d})= {{c 1_A}\over d}$ and 
$\iota_A(a) = {a\over 1}$.
\endrmprofess

\profess{Convention regarding polynomials in several variables}
{\noncommutingvar}
For non-commutative $A$ and $n>1$, $\IntnD(A)$ is a priori not
closed under multiplication and therefore in general not a ring.
With the exception of the following section on $I$-adic continuity,
we will only consider polynomial functions in several variables if the
$D$-algebra $A$ is commutative. 

From section $3$ onward, statements about $\IntnD(A)$ with
unspecified $n$ and $A$ are meant as follows:
if $A$ is commutative, let $n\in\natn_0$, if $A$ is non-commutative, 
assume $n\le 1$.
\endprofess

\rmprofess{Remark}{\continInt}
Note that $K\cap A= D$ implies
$$\IntD(A)\subseteq \Int(D)=\{f\in K[x]\mid f(D)\subseteq D\}$$
and for commutative $A$, also
$$\IntnD(A)\subseteq \Int(D^n)=\{f\in K[\xn]\mid f(D^n)\subseteq D\},$$
and $\IntzeroD(A)= K\cap A=D$.
\endrmprofess

\subhead{Continuity}

$I$-adic continuity of the function $f\colon D^n\rightarrow D$
arising from a classical integer-valued polynomial 
$f\in\Int(D^n)=\{f\in K[\xn]\mid f(D^n)\subseteq D\}$
has been shown, for an arbitrary ideal $I$ of an arbitrary domain 
$D$, in Prop.~1.4 of \cite{CCF00id}. (The proof there is for one 
variable, but clearly generalizable to several variables.)

To establish $I$-adic continuity of integer-valued polynomials
on algebras, we will briefly look at polynomials in several 
non-commuting variables.
If our algebra $A$ is non-commutative, this becomes necessary, even
if we are only interested in integer-valued polynomials in one
variable: if we consider $f(x+y)-f(y)$ as a polynomial in two
variables, and we still want substitution of elements from $A$ 
for $x$ and $y$ to be a homomorphism, we must turn to 
non-commuting variables.

\rmprofess{Definition}{\severalnoncommdef}
Let $D$ be a domain with quotient-field $K$.
Let $K\langle\xn\rangle$ be the free associative $K$-algebra
generated by $\xn$ (in other words, the semigroup-ring $K(S)$,
where $S$ is the free semigroup generated by $\xn$). 
If $A$ is a torsion-free $D$-algebra, we evaluate polynomials
in $K\langle\xn\rangle$ at arguments in $B=K\tensor_D A$ and
thus associate a polynomial function $f\colon B^n\rightarrow B$
to every $f\in K\langle\xn\rangle$. Such polynomials as map
arguments in $A^n$ to values in $A$ we call integer-valued on $A$.
\endrmprofess

From the theory of PID-rings (polynomial identity rings), it is easy
to garner non-trivial examples of integer-valued polynomials in
several non-commuting variables.
For instance, if $p$ is prime and $n\ge 1$, then a polynomial in 
$\ratq\langle x_1,\ldots, x_{np}\rangle$,
but not in $\intz\langle x_1,\ldots, x_{np}\rangle$,
that takes every $np$-tuple of $n\times n$ integer matrices to an 
integer matrix is
$$f(x_1,\ldots,x_{np}) =
\frac{1}{p}\sum_{\pi\in S_{np}} x_{\pi(1)}x_{\pi(2)}\ldots x_{\pi(np)}.$$
(This follows from \cite{Cha88PIp} or \cite{Zal85sAL}.) 
In this paper, we will not consider polynomials in non-commuting 
variables except in the following theorem and its corollaries.

\profess{Theorem}{\noncommcontinuity} 
Let $D$ be a domain with quotient-field $K$ and $A$ a torsion-free
$D$-algebra.
For every $f\in K\langle\xn\rangle$ integer-valued on $A$,
the polynomial function $f\colon A^n\rightarrow A$ is
uniformly $I$-adically continuous for every ideal $I$ of $D$.
\endprofess

\proof
Fix $i$ and let $d$ be the degree of $f$ in $x_i$. We will show that
for every $b\in I^{d}A$ and every $(\an)\in A^n$, 
$$f(a_1,\ldots,a_i+b,\ldots, a_n) - f(a_1,\ldots,a_i,\ldots, a_n)\in IA.
\leqno{(\ast)}$$
For $d=0$, or if $f$ is the zero-polynomial, this is obvious. 
The polynomial
$f$ is uniquely representable as $f=f_1+f_2$, where $x_i$ doesn't occur
in $f_1$, $x_i$ occurs in every monomial in the support of $f_2$, and
$f_2$ has the same degree in $x_i$ as $f$.
Since 
$f_1=f(x_1,\ldots,x_{i-1},0,x_{i+1},\ldots,x_n)$ and $f_2=f-f_1$,
both $f_1$ and $f_2$ are integer-valued. We can show $(\ast)$
separately for $f_1$ and $f_2$. As $(\ast)$ holds for 
$f_1$ and arbitrary $b\in A$, we have reduced to the case $f=f_2$, i.e.,
when $x_i$ occurs in every monomial in the support of $f$.

Also, it suffices to show $(\ast)$ for 
$b=t_d t_{d-1}\ldots t_1 c$ with $t_k\in I$ and $c\in A$, because every
element of $I^{d}A$ is a finite sum of elements of this form.

By considering 
$$g(\xn,z)= f(x_1,\ldots,x_i+z,\ldots, x_n) - f(x_1,\ldots,x_i,\ldots, x_n)$$
we can reduce our task to showing
for every $g(\xn,z)\in K\langle\xn,z\rangle$ of degree $d\ge 1$ in $z$,
which maps $A^{n+1}$ to $A$, satisfies $g(\xn,0)=0$, and such that
$z$ occurs in every monomial in the support of $g$: 
$$g(\an,t_d\ldots t_1c)\in IA\quad \hbox{\rm whenever}\quad
t_1,\ldots,t_d\in I \hbox{\rm\ \ and\ \ }c\in A.\leqno{(\ast\ast)}$$
We show $(\ast\ast)$ by induction on $d$. 
Let $d\ge 1$. Consider
$$h(\xn,y)=g(\xn,t_dy)-t_d^d\,g(\xn,y).$$ 
$h$ satisfies $h(\xn,0)=0$, is of degree at most $d-1$ in $y$, and
$y$ occurs in every monomial in the support of $h$.
Now $h(\an,t_{d-1}\ldots t_1c)\in IA$
for all $t_1,\ldots, t_{d-1}\in I$ and $c\in A$, either by induction
hypothesis or because $h$ is the zero-polynomial. Therefore
$g(\an,t_d t_{d-1}\ldots t_1 c)=
h(\an,t_{d-1}\ldots t_1 c) +t_d^d\,g(\an,t_{d-1}\ldots t_1 c)$ is in $IA$,
for all $t_1,\ldots, t_d\in I$ and $c\in A$.
\endproof

Returning to commuting variables, we conclude:

\profess{Corollary}{\continuity} For any ideal $I$ of $D$, and every 
$f\in\IntnD(A)$ the function $f\colon A^n\rightarrow A$
is uniformly $I$-adically continuous.
\endprofess

\profess{Corollary}{\extendtocompletion}
If $M$ is a maximal ideal of $D$, $\Ahat$  the $M$-adic completion of
$A,$ and $f\in\IntnD(A)$, then the function $f\colon A^n\rightarrow A$ 
defined by $f$ extends uniquely to a $M$-adically continuous function
$f: \Ahat^n\rightarrow \Ahat$.
\endprofess

\subhead{A few technicalities}

This section contains lemmata needed for the investigation of the
spectrum of $\IntD(A)$ and, for commutative $A$, of $\IntnD(A)$. 
From now on, all statements about $\IntnD(A)$ are subject to
the convention: if $A$ is non-commutative, assume $n\le 1$.

\profess{Lemma}{\pseudoprincipal}
Let $D$ be a domain and $P$ a finitely generated prime ideal of height $1$. 
For every non-zero $p\in P$ there exist $m\in \natn$ and $s\in D\takeaway P$ 
such that $sP^m\subseteq pD$.
\endprofess

\proof
Let $p\in P$.
In the localization $D_P$, $P_P$ is the radical of $(q)$ for
every non-zero $q\in P_P$. Therefore,
since $P$ (and hence $P_P$) is finitely generated, there exists
$m\in\natn$ with ${P_P}^m\subseteq pD_P$ and in particular
$P^m\subseteq pD_P$.

The ideal $P^m$ is also finitely generated, by $p_1,\ldots,p_k$, say.
Let $a_i\in D_P$ with $p_i=pa_i$. By considering the fractions
$a_i=r_i/s_i$ (with $r_i\in D$ and $s_i\in D\takeaway P$), and 
setting $s=s_1\cdot\ldots\cdot s_k$, we see that $sP^m\subseteq pD$ 
as desired.
\endproof

\rmprofess{Definition}{\intonsubsetdef}
If $S\subseteq B^n$ and $T\subseteq B$, let
$$\IntnD(S,T)=\{f\in K[\xn]\mid f(S)\subseteq T\}.$$
\endrmprofess

\profess{Lemma}{\IntAPAcontainedinQ}
Let $D$ be a domain and $P$ a finitely generated prime ideal
of height $1$. Then every prime ideal $Q$ of $\IntnD(A)$ with
$Q\cap D=P$ contains $\IntnD(A, PA)$.
\endprofess

\proof
Let $f\in \IntnD(A, PA)$. By dint of Lemma \pseudoprincipal,
there are $m\in\natn$, $p\in P$ and $s\in D\takeaway P$ such
that $sP^m\subseteq pD$.
Then $sf^m\in \IntnD(A, pA)=p\,\IntnD(A)\subseteq Q$. 
As $Q$ is prime and $s\notin Q$, we conclude that $f\in Q$.
\endproof

\profess{Lemma}{\Mafiniteindex}
Let $A$ be a $D$-algebra that is finitely generated as a $D$-module,
$M$ a maximal ideal of finite index in $D$ and $\Ahat$ the $M$-adic
completion of $A$.
For all $a\in \Ahat^n$, the ideal
$\MAhata=\{f\in\IntnD(A)\mid f(a)\in M\Ahat\}$
is of finite index in $\IntnD(A)$.
\endprofess

\proof
Since $A$ is a finitely generated $D$-module, $\Ahat$ is a finitely
generated $\Dhat$-module. As $\Mhat$ is of finite index in $\Dhat$, 
$\Ahat/M\Ahat$ is finite. Let $\IntnD(A)(a)$ denote the image of $a$
under $\IntnD(A)$. Then $\IntnD(A)(a)/(M\Ahat\cap \IntnD(A)(a))$,
as a subring of $\Ahat/M\Ahat$, is finite.

Let $b_1,\ldots,b_m$ be a system of representatives of 
$\IntnD(A)(a)$ modulo $M\Ahat\cap \IntnD(A)(a)$,
and for $1\le i\le m$ let $f_i\in\IntnD(A)$ with $f_i(a)=b_i$.
Then for every $f\in\IntnD(A)$, exactly one of the differences $f-f_i$
is in $\MAhata$, which means that $\{f_1,\ldots,f_m\}$ is a complete
system of residues of $\IntnD(A)(a)$ modulo $\MAhata$.
We have shown 
$[\IntnD(A):\MAhata] = [\IntnD(A)(a):M\Ahat\cap \IntnD(A)(a)]\le
[\Ahat:M\Ahat]$
\endproof

\subhead{Primes lying over a height one maximal ideal of finite index.}

In this section we determine the prime ideals of $\IntD(A)$ (and, for
commutative $A$, of $\IntnD(A)$) lying over a height one prime ideal
of finite index of $D$. 
Prime ideals lying over a prime of infinite index in $D$ will be
characterized in the next section.

\profess{General hypotheses in this section}{\spectrumhypotheses}
\item{\(i\)} $D$ is a domain and $A$ a torsion-free $D$-algebra,
finitely generated as a $D$-module;
\item{\(ii\)} $M$ is a finitely generated maximal ideal
of height $1$ and finite index in $D$;
\item{\(iii\)}
$MA_M\cap A=MA$. (Note that $MA_M\cap A=MA$ is satisfied in two
important cases: if $A$ is a free $D$-module, and if $D\subseteq A$ 
is an extension of Dedekind rings.)
\par
\endprofess

We denote the $M$-adic completions of $D$ and $A$ by $\Dhat$ and
$\Ahat$. Since $M$ is finitely generated, $M$-adic and
$\Mhat$-adic topologies coincide on $\Dhat$ and on $\Ahat$.

\profess{Lemma}{\jeanlucremark} 
The hypotheses of \spectrumhypotheses\ imply
\item{\(iv\)} $M\Ahat\cap A= MA$.
\item{\(v\)} $\Dhat$ and $\Ahat$ are compact.
\endprofess

\proof
\(iv\)
Whenever $M$ is a finitely generated maximal ideal of height $1$ in a 
domain $D$ and $A$ a finitely generated $D$-module, the equality
$M\Ahat\cap A= MA_M\cap A$ holds, by
\cite[Chapter {\rm III}, \S 2.12, Prop.~16]{Bou89CA} combined with
\cite[\S 3.5, Cor.~1 of Prop.~9]{Bou89CA}.

\(v\) The ring $\Dhat$ is compact because $M$ is of finite index and
finitely generated, which implies that all powers of $M$ are
of finite index. $\Ahat$ then is compact because it is finitely 
generated as a $\Dhat$-module.
\endproof

\rmprofess{Notation}{\specnotation}
Recall our convention that we only allow $n>1$ in $\IntnD(A)$ 
if $A$ is commutative; for non-commutative $A$, $n=1$ is assumed.

The image of $a\in\Ahat$ under $\IntnD(A)$ and $\IntD(A)$ we denote
as follows:
$$\IntnD(A)(a)=\{f(a)\mid f\in\IntnD(A)\} \quad\hbox{\rm and}\quad
\IntD(A)(a)=\{f(a)\mid f\in\IntD(A)\}.$$

If $M$ is a maximal ideal of $D$ and $a\in \Ahat^n$ let 
$$\MAhata=\{f\in \IntnD(A)\mid f(a)\in M\Ahat\}.$$

If $P$ is an ideal of a commutative ring between $\IntnD(a)$ and $\Ahat$, 
let
$$\Pa=\{f\in\IntnD(A)\mid f(a)\in P\}.$$
\endrmprofess

\comment{
Since $M\Ahat\cap A=MA$, $\MAhata$ with $a\in A$ becomes
$$\MAa=\{f\in \IntnD(A)\mid f(a)\in MA\}.$$
} 

\profess{Lemma}{\QcontainsMa}
Under the hypotheses of \spectrumhypotheses, let $Q$ be a prime ideal
of $\IntnD(A)$ lying over $M$. 
Then there exists $a\in\Ahat$ such that $$\MAhata\subseteq Q.$$
In particular, $Q$ is of finite index.
\endprofess

\proof
Suppose $Q$ does not contain any $\MAhata$. Then for every $a\in\Ahat^n$
there exists $f\in \MAhata\takeaway Q$. By Theorem \continuity, $f$ is
$M$-adically continuous, so there exists an $M$-adic neighborhood $U$ of
$a$, such that $f(U)\subseteq M\Ahat$. By compactness of $\Ahat^n$, there
exist finitely many $a_i$ such that the corresponding neighborhoods
$U_{i}$ cover $\Ahat^n$. Let $g$ be the product of the polynomials
$f_i\in (M\Ahat)_{a_i}\takeaway Q$ with $f_i(U_i)\subseteq M\Ahat$.
For all $a\in A$, $g(a)\in M\Ahat\cap A= MA_M\cap A=MA$ 
(by Lemma \jeanlucremark).
Since $Q$ is prime, $g\notin Q$, and yet $g\in  \IntnD(A, MA)$, a 
contradiction to Lemma \IntAPAcontainedinQ. 
We have shown that $Q$ contains some $\Mhata$, which is of finite
index by Lemma \Mafiniteindex.
\endproof

In the special case $A=D$, the preceding lemma already concludes
the characterization of primes of $\IntnD(D)$ lying above
a maximal ideal $M$ of finite index in $D$ (a result of Chabert 
\cite{Cha77IP}), because then $\MAhata$ is
$\(M\Dhat)_a=\{f\in\IntnD(D)\mid f(a)\in\Mhat\}$, a prime ideal of
finite index, and hence maximal.
Chabert, however, showed the other inclusion, $Q\subseteq \MAhata$,
after first showing independently that $Q$ must be maximal
\cite[Prop.~V.2.2]{CaCh97ivp}.

\profess{Lemma}{\QMaP} 
Under the hypotheses of \spectrumhypotheses,
let $Q$ be a prime ideal of $\IntnD(A)$ lying over $M$. 
For every $a\in\Ahat$ such that $\MAhata\subseteq Q$,
there exists a maximal ideal $P$ of $\IntnD(A)(a)$
such that $Q=\Pa$.
\endprofess

\proof
Since $\Ahat/ M\Ahat$ is a finite ring, 
$\IntnD(A)(a)/(\IntnD(A)(a)\cap M\Ahat)$ is
a finite commutative ring. Let $P_1,\ldots, P_k$ be the 
maximal ideals of $\IntnD(A)(a)$ containing $\IntnD(A)(a)\cap M\Ahat$. 

Suppose $Q$ is not contained in any $(P_i)_a$, for $1\le i\le k$. 
Then, by prime avoidance, $Q(a)=\{f(a)\mid f\in Q\}$ is not contained
in $\bigcup_{i=1}^k P_i$. 
Let $f\in Q$ such that $f(a)$ is not in any $P_i$. 
Then the residue class of $f(a)$ is a unit in
$\IntnD(A)(a)/(\IntnD(A)(a)\cap M\Ahat)$.

Replacing $f$ by a suitable power of $f$
(using the fact that the group of units of 
$\IntnD(A)(a)/(\IntnD(A)(a)\cap M\Ahat)$ is finite) 
we see that there exists $f\in Q$ with $f(a)\congr 1$ mod $M\Ahat$.
It follows that $1-f\in \MAhata\subseteq Q$ and therefore $1\in Q$,
a contradiction.
\endproof

\profess{Theorem}{\maxideals}
Let $D$ be a domain, $A$ a torsion-free $D$-algebra finitely
generated as a $D$-module, $M$ a finitely generated
maximal ideal of $D$ of finite index and height one,
such that $MA_M\cap A=MA$.

The prime ideals of $\IntnD(A)$ lying over $M$ are
precisely the ideals of the form 
$$\Pa=\{f\in\IntnD(A)\mid f(a)\in P\},$$
where $a\in\Ahat$ (the $M$-adic completion of $A$), and $P$ is a maximal
ideal of $\IntnD(A)(a)$ \(the image of $a$ under $\IntnD(A)$\) with
$P\cap D=M$. In particular, all primes of $\IntnD(A)$ lying over $M$
are of finite index.
\endprofess

\proof
There exist primes of $\IntnD(A)$ lying over $M$, because 
$\MAhata$ with $a\in \Ahat$ is a proper ideal of $\IntnD(A)$
containing $M$.

If $Q$ is a prime ideal of $\IntnD(A)$ lying over $M$, then
Lemma \QcontainsMa\ shows that there exists an element $a\in\Ahat$
such that $\Mhata$ is contained in $Q$, and that $Q$ is of finite index. 
It then follows from Lemma \QMaP\ that $Q=P_a$ for some maximal ideal $P$
of $\IntnD(A)(a)$ satisfying 
$M\subseteq \IntnD(A)(a)\cap M\Ahat\subseteq P$,
and hence $P\cap D=M$.

Conversely, if $P$ is a maximal ideal of $\IntnD(A)(a)$, 
then $\IntnD(A)/P_a$ is isomorphic to $\IntnD(A)(a)/P$.
Therefore $P_a$ is a maximal ideal of $\IntnD(A)(a)$.
\endproof

\comment{
\profess{Lemma}{\finiteringext}
Let $R\subseteq S$ be finite commutative rings. Then every maximal
ideal $Q$ of $R$ lies under a maximal ideal $P$ of $S$.
\endprofess
 
\proof
Let $Q$ be a maximal ideal of $R$. Then $Q$ does not contain any
unit $u$ of $S$, because otherwise $Q$ would contain $1_S=1_R$, which
is a power of $u$. Therefore $Q$ is a ring (possibly without identity)
contained in the union of all maximal ideals of $S$. By prime avoidance,
$Q$ is contained in some maximal ideal $P$ of $S$. Now $QS\subseteq P$,
and, by maximality of $Q$, we see that $P\cap R=Q$.
\endproof
}

It may happen that we do not know the exact image of $a\in\Ahat$
under $\IntnD(A)$, but do know a commutative ring $R_a$ between
$\IntnD(A)(a)$ and $\Ahat$.
In this case we should remember that 
$R_a/(R_a\cap M\Ahat)$ is a subring of the finite ring
$(\Ahat/ M\Ahat)$, and that therefore 
$\IntnD(A)(a)/(\IntnD(A)(a)\cap M\Ahat)\subseteq (R_a/R_a\cap M\Ahat)$
is an extension of finite commutative rings. Since extensions
of finite commutative rings satisfy ``lying over'', every prime ideal
of $\IntnD(A)(a)$ comes from a prime ideal of $R_a$, and we conclude:

\profess{Corollary}{\commoverring} 
Under the hypotheses of of Theorem \maxideals, suppose we have, for
every $a\in\Ahat$, a commutative ring $R_a$ with
$\IntnD(A)(a)\subseteq R_a\subseteq \Ahat$.

Then the prime ideals of $\IntnD(A)$ are precisely the ideals of
the form $$P_a=\{f\in \IntnD(A)\mid f(a)\in P\},$$
where $a\in\Ahat$ and $P$ is a maximal ideal of $R_a$ lying over
$M$.
\endprofess

\comment{
There is a fine point to observe in \maxideals: if a
maximal ideal of $\Intn(A)$ is of the form $Q=P_a$, where  
$a\in\Ahat$ and $P$ is a maximal ideal of $\Imn(a)$ lying over $M$,
it does not automatically follow that $Q$ contains $M_a$.
It may happen that $Q=M_{a,P}=M_{b,P'}$, with $M_b\subseteq Q$, but
$M_a\not\subseteq Q$. 
}

If $A$ is a commutative $D$-algebra, we can take $R_a=\Ahat$ in 
Corollary \commoverring\ for all $a\in\Ahat$, and we get the
following simpler characterization of the primes of $\IntnD(A)$
lying over $M$:

\profess{Theorem}{\commAmaxideals}
Let $D$ be a domain, $A$ a commutative torsion-free $D$-algebra
finitely generated as a $D$-module, $M$ a finitely generated
maximal ideal of $D$ of finite index and height one,
such that $MA_M\cap A=MA$.

Then every prime ideal of $\IntnD(A)$ lying over $M$ is of the form 
$$P_a=\{f\in \IntnD(A)\mid f(a)\in P\},$$
for some $a\in\Ahat$ (the $M$-adic completion of $A$)
and $P$ a maximal ideal of $\Ahat$ lying over $M$. In particular,
every prime of $\IntnD(A)$ lying over $M$ is of finite index.
\endprofess

\subhead{Primes lying over prime ideals of infinite index}

\profess{Lemma}{\infindexprime}
Let $A$ be a torsion-free $D$-algebra with $K\cap A=D$.
Let $P$ be a prime ideal of infinite index in $D$ and $n\in\natn$.
Then $\IntnD(A)\subseteq D_P[\xn]$.
\endprofess

\proof
More generally, we show that for any prime ideal $P$, a polynomial 
$f\in\IntnD(A)$ of degree less than $[D : P]$ in every 
individual variable is in $D_P[\xn]$. We use induction on $n$.
The case $n=0$ is trivial: $\Int^{0}_D(D)= K\cap A=D\subseteq D_P$.

For $n>0$ consider $f$ as a
polynomial in $x_n$ with coefficients in $K[x_1,\ldots,x_{n-1}]$.
Let $s\le [D : P]$ such that $f$ is of degree strictly less 
than $s$ in each $x_j$.

Choose $d_1,\ldots,d_s\in D\subseteq A$ pairwise incongruent mod $P$.
For every $i$, the value of $f$ at $d_i$ (substituted for $x_n$) is
a polynomial in $\Int^{n-1}_D(A)$ of degree less than $s$ in each
variable. Therefore 
$f(x_1,\ldots,x_{n-1},d_i)\in D_P[x_1,\ldots,x_{n-1}]$, by induction
hypothesis.

Let $g\in K(x_1,\ldots,x_{n-1})[x_n]$ be the Lagrange interpolation
polynomial with $g(d_i)=f(x_1,\ldots,x_{n-1},d_i)$ ($1\le i\le s$);
then $g\in D_P[x_1,\ldots,x_{n-1}][x_n]$. 
Since a polynomial (with coefficients in a domain) of degree less than
$s$ is determined by its values at $s$ different arguments, we must
have $f=g\in D_P[x_1,\ldots,x_{n-1}][x_n]$.
\endproof

\profess{Corollary}{\infiniteindextrivial}
Let $A$ be a torsion-free $D$-algebra with $K\cap A=D$.
If all maximal ideals of $D$ are of infinite index, then
$\IntnD(A)=D[\xn]$.
\endprofess

Alternatively, we could have deduced the previous lemma from the
corresponding fact for the ring of integer-valued polynomials 
over $D$ 
\cite[Prop.~I.3.4, XI.1.10]{CaCh97ivp},
since after all 
$\IntnD(A)\subseteq \Int(D^n)=\{f\in K[\xn]\mid f(D^n)\subseteq D\}$.
\comment{
The corollary raises the question when $\IntnD(A)\ne D[\xn]$. For
Noetherian $D$, a criterion for $\IntnD(A)\ne D[\xn]$, similar to
the criterion for $\Int(D)\ne D[x]$ (\cite{CaCh97ivp} Thm I.3.14)
is given in \cite{Fri--mrs}.
}

\profess{Lemma}{\primesoverinfiniteindex}
Let $A$ be a torsion-free $D$-algebra with $K\cap A=D$.
Let $P$ be a prime ideal of infinite index in $D$ and $n\in\natn$.
Then the prime ideals of $\IntnD(A)$ lying over $P$ are precisely
those of the form $Q\cap \IntnD(A)$, where $Q$ is a prime ideal of
$D_P[\xn]$ containing $PD_P[\xn]$.
\comment{
In particular, the prime ideals of $\IntD(A)$ lying
over $(0)\ideal D$ are precisely those of the form 
$f(x)K[x]\cap \Int(A)$ where $f$ is an irreducible polynomial in
$K[x]$ and for non-zero $P\ideal D$ of infinite index the prime ideals
of $\Int(A)$ lying over $P$ are precisely those of the form 
$(f(x)D_P[x]+PD_P[x])\cap \IntD(A)$ where $f$ represents an
irreducible polynomial in $(D_P/PD_P)[x]$.
}
\endprofess

\proof
As $D[\xn]\subseteq \IntnD(A)\subseteq D_P[\xn]=D[\xn]_{(D\setminus P)}$,
we have $$D_P[\xn]=\IntnD(A)_{(D\setminus P)}$$ 
and therefore a bijective correspondence (given by lying over) exists
between prime ideals of $\IntnD(A)$ whose intersection with $D$ is
contained in $P$ and prime ideals of $D_P[\xn]$. 
\endproof

\profess{Theorem}{\dimension}
Let $D$ be a Noetherian one-dimensional domain with finite
residue fields and $A$ a torsion-free $D$ algebra, finitely
generated as a $D$-module, such that for every maximal ideal
$M$, $MA_M\cap D=M$. If $A$ is commutative, let $n\in\natn$,
for non-commutative $A$ restrict to $n=1$. Then $\IntnD(A)$
is $(n+1)$-dimensional.
\endprofess

\proof
By Lemma \primesoverinfiniteindex,
the prime ideals of $\IntnD(A)$ lying over $(0)$ all come from
prime ideals of $K[\xn]$.
The primes of $\IntnD(A)$ lying over a maximal ideal
$M$ are all maximal and hence mutually incomparable.
So $\dim(\IntnD(A))\le n+1$.
If $M$ is a maximal ideal of $D$ and $d=(d_1,\ldots,d_n)\in D^n$,
let ${\cal Q}_k$ be the ideal of $K[\xn]$ generated by 
$(x_1-d_1),\ldots, (x_k-d_k)$. 
Then a chain of primes of length $n+1$ of $\IntnD(A)$ is given by
$Q_0=(0)$, $Q_k={\cal Q}_k\cap \IntnD(A)$, for $k=1,\ldots,n$ and 
$Q_{n+1}=P_d=\{f\in\IntnD(A)\mid f(d_1,\ldots,d_n)\in P\}$, where
$P$ is a maximal ideal of the image of $d$ under $\IntnD(A)$.
\endproof

\comment{
Recall that a d-ring is a domain $D$ such that every non-constant
polynomial in $D[x]$ has a zero modulo some maximal ideal of $D$
\cite{CaCh97ivp}, Prop.~VII.2.3.
Every domain of characteristic $0$, finitely generated as $\intz$-algebra,
is a d-ring, \cite{CaCh97ivp}, Prop.~VII.2.9. 
}

\rmprofess{Remark}{\nomaxunitary}
If $D$ is a Noetherian domain of characteristic $0$, finitely generated
as a $\intz$-algebra, such as, for instance, the ring of integers $\OK$
in a number field, then no maximal ideal of $\IntnD(A)$ lies over $(0)$
of $D$.  This holds for $A=D$ by 
\cite[Prop.~XI.3.4.]{CaCh97ivp}, 
and carries over to $\IntnD(A)$, since $\IntnD(A)\subseteq \Int(D^n)$.
Every prime ideal $P$ of $\IntnD(A)$ coming from a prime ideal of
$K[\xn]$ is contained in a maximal ideal $Q$ of $\Int(D^n)$ lying
over a maximal ideal $M$ of $D$, and $Q\cap \IntnD(A)$ then
properly contains $P$.
\endrmprofess

\rmprofess{Remark}{\lengthofchains}
Even if there are no maximal ideals lying over $(0)$, maximal chains
of primes of $\IntnD(A)$ are not necessarily of length $n+1$. For 
instance, a maximal ideal of $\IntnD(A)$ may well have height $1$ if 
it is of the form $\MAhata$ for $a=(a_1,\ldots,a_n)\in\Ahat^n$
with $a_1,\ldots,a_n$ algebraically independent over $K$. 
\endrmprofess
 
\subhead{Integer-valued polynomials on matrix algebras}

Theorem \maxideals\ characterizes the spectrum of the ring $\IntD(A)$,
provided we know the images of elements of $M$-adic completions
of $A$ under $\IntD(A)$. We will now determine these images in
the case $A=M_n(D)$. Note that all the technical hypotheses in
this section are certainly satisfied if $D=\OK$ is the ring of 
integers in a number field.

\profess{Fact}{\whenMint} {\rm \cite[Lemma 2.2]{Fri05chh}}
Let $D$ be a domain and $f(x)=g(x)/d$ with $g\in D[x]$,
$d\in D\setminus \{0\}$. Then $f\in\IntD(M_n(D))$ if and only if $g$ is
divisible modulo $dD[x]$ by all monic polynomials in $D[x]$ of
degree $n$.
\endprofess

\profess{Proposition*}{\simplerwhenMint}
Let $D$ be a domain
such that the intersection of the maximal ideals of $D$ of finite index 
is $(0)$, and $f(x)=g(x)/d$ with
$g\in D[x]$, $d\in D\setminus \{0\}$. Then $f\in\IntD(M_n(D))$ if and
only if $g$ is divisible modulo $dD[x]$ by all monic irreducible 
polynomials in $D[x]$ of degree $n$.
\endprofess

\proof
In view of Fact \whenMint, it suffices to show for every 
$d\in D\setminus\{0\}$ and $h\in D[x]$ monic of degree $n$,
that there exists $k\in D[x]$ monic of degree $n$, irreducible in
$D[x]$ and congruent to $h$ mod $dD[x]$. 
We may choose a maximal ideal $P$ of finite index with $d\notin P$,
and use Chinese remainder theorem on the coefficients of $h$ to find 
$k\in D[x]$, monic of degree $n$, congruent to $h$ mod $dD[x]$ and 
irreducible in $(D/P)[x]$. 
\endproof

We are now able to identify a polynomially dense subset of $M_n(D)$
consisting of companion matrices. They are often easier to work with
than general matrices, because their characteristic polynomial is
also their minimal polynomial. 

\profess{Theorem*}{\compdense}
Let ${\cal C}_n$ be the set of companion matrices of monic polynomials
of degree $n$ in $D[x]$ and ${\cal I}_n\subseteq {\cal C}_n$ the subset
of companion matrices of irreducible polynomials.
If $D$ is any domain, 
$$\IntD(M_n(D))= \IntD({\cal C}_n,M_n(D)).$$
If $D$ is a domain
such that the intersection of the maximal ideals of
finite index is $(0)$, such as, for instance,
a Dedekind domain with infinitely many maximal ideals of finite index,
then
$$\IntD(M_n(D))= \IntD({\cal I}_n,M_n(D)).$$
\endprofess

\proof
Let $f\in \IntD({\cal C}_n,M_n(D))$, $f(x)=g(x)/d$ with $g\in D[x]$,
$d\in D$.
Since $g$ maps every $C\in {\cal C}_n$ to $M_n(dD)$, $g$ is divisible
mod $dD[x]$ by every monic polynomial in $D[x]$ of degree $n$. (This is
so because $f$ is still the minimal polynomial of its companion matrix
when everything is viewed in $D/dD$.) By Fact~\whenMint,
$f\in \IntD(M_n(D))$. This shows $\IntD({\cal C}_n,M_n(D))\subseteq
\IntD(M_n(D))$. The reverse inclusion is trivial. The argument for
${\cal I}_n$ is similar, using Prop.~\simplerwhenMint.
\endproof

\profess{Theorem}{\mimage}
Let $D$ be a domain and $C\in M_n(D)$. Let
$$\Im(C)=\{f(C)\mid f\in\IntD(M_n(D))\}
\quad{\rm\ and\ }\quad
D[C]=\{f(C)\mid f\in D[x]\}.$$ 
Then $\Im(C)=D[C]$.
\endprofess

\proof
Consider $f\in\IntD(M_n(D))$; $f(x)=g(x)/d$ with $g\in D[x]$ and
$d\in D\setminus\{0\}$. We know that $g$ is divisible modulo $dD[x]$
by every monic polynomial in $D[x]$ of degree $n$. Dividing $g$ by
$\chiC$, the characteristic polynomial of $C$, we get
$$g(x)=q(x)\chiC(x)+dr(x)$$ with $q,r\in D[x]$ and we see that
that $f(C)=r(C)$. Thus $\Im(C)\subseteq D[C]$. The reverse inclusion
is clear, since $D[x]\subseteq\IntD(M_n(D))$.
\endproof

\rmprofess{Definition}{\defanairred}
A local domain is called {\em analytically irreducible} if its completion 
is also a domain.
\endrmprofess

\profess{Lemma}{\anairred}
Let $M$ be a maximal ideal of finite index in a domain $D$.
Then the following are equivalent
\item{\rm (1)}
$\bigcap_{n=1}^\infty M^n = (0)$ 
and $D_M$ is analytically irreducible
\item{\rm (2)}
for every non-zero $d\in D$,
cancellation of $d$ is uniformly $M$-adically continuous.
\endprofess

\proof
$(1\Rightarrow 2)$
We have to show: for every non-zero $d\in D$, for every $m\in\natn$ there
exists $k\in\natn$ such that for all $c\in D$: $dc\in M^{k}$ implies
$c\in M^m$.
Indirectly, suppose there exists $d\in D\setminus \{0\}$ and $m\in\natn$,
such that for every $k\in\natn$ there is some $c_k\in D$ with
$dc_k\in M^{k}$ and $c_k\notin M^m$.
Since $\Dhat$ is compact and satisfies first countability axiom,
$(c_k)$ has a convergent subsequence. Its limit $c\in\Dhat$ satisfies
$c\notin M^m$, and hence $c\ne 0$, and also for all $k$, $dc\in M^k$,
which implies $dc=0$. We have shown the existence of zero-divisors in
$\Dhat$. 

$(2\Rightarrow 1)$ is easy.
\endproof

\profess{Theorem}{\complimage}
Let $D$ be a domain, $M$ a maximal ideal of finite index of $D$
such that $\bigcap_{n=1}^\infty M^n = (0)$ and $D_M$ is analytically
irreducible. Let $\Dhat$ be the $M$-adic completion of $D$,
$C\in M_n(\Dhat)$, and $\IntD(M_n(D))(C)$ the image of $C$
under $\IntD(M_n(D))$. Then $$\IntD(M_n(D))(C)\subseteq \Dhat[C].$$
\endprofess 

\proof Let $f\in \IntD(M_n(D))$, $f(x)=g(x)/d$ with $g\in D[x]$, $d\in D$.
For every $m\in\natn$ let $k_m\in \natn$ such that for all $c\in D$,
$cd\in M^{k_m}$ implies $c\in M^m$.
Let $E_1=(e_{ij}^{(1)})$ and $E_2=(e_{ij}^{(2)})$ be matrices in $M_n(D)$
with characteristic polynomials $\chi_1$ and $\chi_2$. Then
$g(x)=q_i(x)\chi_i(x) + dr_i(x)$ with $q_i, r_i\in D[x]$ for $i=1,2$.
If $e_{ij}^{(1)}\congr e_{ij}^{(2)}$ mod $M^{k_m}$, then
$dr_1\congr dr_2$ mod $M^{k_m}D[x]$, and therefore $r_1\congr r_2$ mod 
$M^mD[x]$. We can therefore $M$-adically approximate $C$ by matrices
$C_i$ with $f(C_i)=s_i(C_i)$ with $s_i\in D[x]$, $\deg s_i<n$, such
that the $s_i$ converge towards a polynomial $s\in \Dhat[x]$ with
$\deg s<n$ and $f(C)=s(C)$.
\endproof

\profess{Corollary}{\Dedecomplimage}
Let $D$ be a Dedekind domain, $M$ a maximal ideal of finite index,
$\Dhat$ the $M$-adic completion of $D$, and $C\in M_n(\Dhat)$. Then
$$\IntD(M_n(D))(C)\subseteq \Dhat[C].$$
\endprofess

\comment{
The maximal ideals of $\Dhat[x]$ are known, of course. They are
the ideals generated by $\Mhat$ and a polynomial in $\Dhat[x]$ that
is irreducible in $(\Dhat/\Mhat)[x]$. 
Now if $C\in M_n(\Dhat)$ and
$N=\{f\in \Dhat[x]\mid f(C)=0\}$ is the null ideal of $C$ in
$\Dhat[x]$, then substitution of $C$ for $x$ gives an isomorphism 
$\Dhat[x]/N\simeq \Dhat[C]$, and the maximal ideals of $\Dhat[C]$ all
come from maximal ideals of $\Dhat[x]$ containing $N$. 

If we want $N$ to be a principal ideal of $\Dhat[x]$, for every
matrix $C\in M_n(\Dhat)$, for every $n$, we must limit ourselves
to integrally closed $\Dhat$ \cite{Fri04ICMP}.

If $\Dhat$ is an integrally-closed domain, then $N$ is a principal
ideal of $\Dhat[x]$, generated by a monic polynomial $f\in \Dhat[x]$,
the minimal polynomial of $C$.
Then $\Dhat[x]/(f)\simeq \Dhat[C]$, and the maximal ideals of 
$\Dhat[x]/(f)$ are exactly the ideals generated by $M$ and a
modulo $\Mhat[x]$ irreducible factor of $f$.
}

\subhead{Integer-valued polynomials with matrix coefficients}

%
While we have been 
investigating the ring $\IntD(M_n(D))$ of polynomials in $K[x]$
mapping matrices in $M_n(D)$ to matrices in $M_n(D)$, 
Werner \cite{Wer12IVM} has been studying the set, let's call it
$\IntD[M_n(D)]$ with square brackets, of polynomials with coefficients
in the non-commutative ring $M_n(K)$ mapping matrices in $M_n(D)$ to 
matrices in $M_n(D)$. Without substitution homomorphism, it is not
a priori clear that this set is closed under multiplication, but
Werner \cite{Wer12IVM} has shown that it is, and so $\IntD[M_n(D)]$
is actually a ring between $(\MnD)[x]$ and $(\MnK)[x]$.

Also in \cite{Wer12IVM}, Werner proves that every ideal of $\MMint$
can be generated by elements of $K[x]$. 
Using the idea of his proof, one can show more:
$\MMint$ is isomorphic to the algebra of $n\times n$ matrices over
$\Mint$. Since every prime ideal of a matrix ring is just the set of
matrices with entries in a prime ideal of the ring, we get a 
description of the spectrum of $\MMint$ as a byproduct of our
description of the spectrum of $\Mint$ 
in the previous section. We recall the definition of prime ideal for
non-commutative rings:

\rmprofess{Definition}{\primeidealdef}
We call a two-sided ideal $P\ne R$ of a (not necessarily commutative) 
ring with identity $R$ a prime ideal, if,
for all ideals $A$, $B$ of $R$,
$$AB\subseteq P\Longrightarrow A\subseteq P {\rm\ \ or\ \ } B\subseteq P,$$
or equivalently, if, for all $a,b\in R$
$$aRb\subseteq P \Longrightarrow a\in P {\rm\ \ or\ \ } b\in P.$$
For commutative $R$ this is equivalent to the (in general stronger)
condition: for all $a,b\in R$,
$$ab\in P \Longrightarrow a\in P {\rm\ \ or\ \ } b\in P.$$
\endrmprofess

\profess{Theorem}{\matrixivpoly} 
Let $D$ be a domain with quotient field $K$, and 
$$ \eqalign{
\Mint=\{f\in K[x]\mid \forall C\in M_n(D):\, f(C)\in  M_n(D)\}\cr
\MMint=\{f\in (M_n(K))[x] \mid \forall C\in M_n(D):\, f(C)\in  M_n(D)\}.\cr
}$$
We identify $\MMint\subseteq (M_n(K))[x]$ with its isomorphic image
in $M_n(K[x])$ under
$$
\varphi\colon (M_n(K))[x]\rightarrow M_n(K[x]),\quad
\sum_k (a_{ij}^{(k)})_{\scriptscriptstyle1\le i,j\le n}\, x^k \mapsto 
\left( \sum_k a_{ij}^{(k)} x^k \right)_{\scriptscriptstyle\!\!1\le i,j\le n}
$$
Then $\MMint = M_n(\Mint)$.
\endprofess

\proof
Note that $K[x]$ is embedded in $M_n(K[x])$ as the subring of scalar
matrices $g(x)I_n$, and in $M_n(K)[x]$ as the subring of polynomials
$g(x)$ whose coefficients are scalar matrices $rI_n$, with $r\in K$.
Clearly, $\MMint\cap K[x] =\Mint$.

Let $C=(c_{ij}(x))\in \MMint\subseteq M_n(K[x])$.
Let $e_{ij}$ be the matrix in $M_n(D)$ with $1$ in position $(i,j)$
and zeros elsewhere; then 
$e_{ij}Ce_{kl}$ has $c_{jk}(x)$ in position $(i,l)$ and 
zeros elsewhere. Also, $e_{ij}Ce_{kl}\in \MMint$, since $\MMint$
is a ring containing $\MnD$.
So $\sum_{i=1}^n e_{ij}Ce_{ki}=c_{jk}(x)I_n\in \MMint$.
Therefore $c_{jk}(x)\in K[x]\cap \IntD[M_n(D)] = \IntD(M_n(D))$ for
all $(j,k)$, and hence $\MMint\subseteq M_n(\IntD(M_n(D)))$.

Conversely, if $f\in\IntD(M_n(D))$ then $f(x)I_n\in \MMint$. Therefore,
$e_{ik}f(x)I_ne_{kl}$, the matrix containing $f(x)$ in position
$(i,l)$ and zeros elsewhere, is in $\MMint$, for arbitrary $(i,l)$.
By summing matrices of this kind we see that
$M_n(\IntD(M_n(D)))\subseteq \MMint$.
\endproof

For any ring $R$ with identity the ideals of $R$ are in bijective 
correspondence with the ideals of $M_n(R)$ by $I\mapsto M_n(I)$,
and restriction to prime ideals gives a bijection between the
spectrum of $R$ and the spectrum of $M_n(R)$. So we conclude:

\profess{Corollary}{\matrixringspec}
Let $\IntD(M_n(D))$ and $\IntD[M_n(D)]$ as in
the preceding theorem. Under the identification of $\IntD[M_n(D)]$
with its isomorphic image in $M_n(K[x])$,
\item{\rm (1)}
The two-sided ideals of $\IntD[M_n(D)]$ are precisely the sets of
the form $M_n(I)$, where $I$ is an ideal of $\IntD(M_n(D))$.
\item{\rm (2)}
The two-sided prime ideals of $\IntD[M_n(D)]$ are precisely the
sets of the form $M_n(P)$, where $P$ is a prime ideal of
$\IntD(M_n(D))$.
\endprofess

\medskip

{\bn Acknowledgments.}
The author wishes to thank J.-L.~Chabert for helpful hints, in particular 
Lemma \jeanlucremark\kern2pt \(iv\), and P.-J.~Cahen for improving
Lemma~\pseudoprincipal.
\medskip

{\bn * Corrigenda.} The published version of this paper in 
J.~Algebra 373 (2013) 414-425 unfortunately contains
a typo and an error.

The typo: in the last sentence of the introductory first 
paragraph of section 7, instead of
``a ring between $\IntD(\MnD)$ and $\IntD(\MnK)$''
(as in the paper published in J Algebra) it should be 
``a ring between $(\MnD)[x]$ and  $(\MnK)[x]$''.
This has been corrected in this post-preprint.

The error: in Prop.~6.2 and Thm.~6.3, instead of
``a domain with zero Jacobson radical'' (as in the
paper published in J Algebra) we need the stronger
assumption ``a domain such that the intersection of
all maximal ideals of finite index is trivial''.
This, too, has been corrected in this post-preprint.
\closeout\xrefs

\vskip-5ex
\nonumsubhead{References}%
\vskip-1ex
\bibliography{mint}
\bibliographystyle{siamese}

\bigskip\bigskip
{ \obeylines\frenchspacing
Institut f\"ur Mathematik A
Technische Universit\"at Graz
Steyrergasse 30
A-8010 Graz, Austria
\smallskip
\tt
frisch@TUGraz.at
} 

\bye